\documentclass{amsart}
\newcommand{\cT}{{\mathcal T}}
\newcommand{\cH}{{\mathcal H}}
\newcommand{\cN}{{\mathcal N}}
\newcommand{\fg}{{\mathfrak g}}
\newcommand{\uc}{\underline{A}}
\newcommand{\codim}{\mbox{codim}}
\newcommand{\uH}{\underline{H}}
\newcommand{\uN}{\underline{N}}
\newcommand{\tuH}{\tilde{\underline{H}}}
\newcommand{\BZ}{\mathbb{Z}}
\newcommand{\BR}{\mathbb{R}}
\newcommand{\BQ}{\mathbb{Q}}
\newcommand{\bN}{\mathbb{N}}
\newcommand{\ox}{\overline{x}}
\newcommand{\oy}{\overline{y}}

\newcommand{\End}{\mbox{End}}
\newcommand{\Ext}{\mbox{Ext}}
\newcommand{\Ho}{\mbox{H}}

\title{Dimensions of quantized tilting modules}
\author{Viktor Ostrik}
\address{Department of Mathematics, Massachusetts Institute of Technology,
Cambridge MA, 02139}
\date{November 2000}
\email{ostrik@math.mit.edu}

\begin{document}

\begin{abstract}
Let $U$ be the quantum group with divided powers at $p-$th root of unity for
prime $p$. For any two-sided cell $A$ in the corresponding affine Weyl group
one associates tensor ideal in the category of tilting modules over $U$. In 
this note we show that for any cell $A$ there exists tilting module $T$
from the corresponding tensor ideal such that biggest power of $p$ which
divides $\dim T$ is $p^{a(A)}$ where $a(A)$ is Lusztig's $a-$function. This
result is motivated by the Conjecture of J.~Humphreys in \cite{H}.
\end{abstract}

\maketitle

\section{Introduction}

Let $G$ be a simply connected algebraic group. In \cite{LuC4} G.~Lusztig
proved existense of bijection between two finite sets: the set of two-sided
cells in the affine Weyl group attached to $G$ (this set is defined 
combinatorially) and the set of unipotent $G-$orbits. The proof in \cite{LuC4}
is quite involved and this bijection remains rather mysterious. In \cite{H}
J.~Humphreys suggested a natural conjectural construction for Lusztig's
bijection using cohomology of tilting modules over algebraic groups in 
characteristic $p>0$ or, similarly, over quantum groups at a root of unity.
In \cite{O1, O2} the author proved some partial results towards this Humphreys'
Conjecture. Now this Conjecture is known to be true in a quantum group case
thanks to the (unpublished) work of R.~Bezrukavnikov. 

For an element $w$ of a finite Weyl group $W_f$ one defines number
$a(w)$ as Gelfand-Kirillov dimension of highest weight simple module 
$L(w\cdot 0)$ over the corresponding semisimple Lie algebra. Generalizing
this G.~Lusztig defined $a-$function on any Coxeter group, see \cite{LuC1}.
The $a-$function takes constant value on any two-sided cell and appears to
be very useful for the theory of cells in Coxeter groups.
J.~Humphreys suggested that his construction of Lusztig's bijection is 
compatible with the theory of $a-$function in the following way: dimension
of any tilting module corresponding to two-sided cell $A$ is divisible by
$p^{a(A)}$ and generically is not divisible by higher power of $p$ (here
$p$ is characteristic of field in a case of algebraic group and order of
root of unity in a case of quantum group). In \cite{O2} the author proved that
the first statement (divisibility of dimensions by $p^{a(A)}$) is a consequence
of Humphreys' Conjecture, so it is a consequence of Bezrukavnikov's work. 
The second statement (generical indivisibility by $p^{a(A)+1}$) seems to be
harder. The main result of this note is that in quantum group case for any 
cell $A$ there exists
tilting module $T$ corresponding to $A$ such that its dimension is not
divisible by $p^{a(A)+1}$. So we determine $p-$component of dimension of
certain tilting modules what seems to be of some interest independently of
Humphreys' Conjecture.  

We will follow the notations of \cite{6}. 
Let $(Y,X,\ldots)$ be a simply connected root datum of finite type. Let $p$
be a prime number bigger than the Coxeter number $h$. 
Let $\zeta$ be a primitive $p$-th root of unity in ${\mathbb C}$. 
Let $U$ be the quantum group with divided
powers associated to these data. Let
$\cT$ be the category of {\em tilting} modules over $U$, see e.g.
\cite{An}. Recall that any tilting module is a sum of indecomposable
ones, and indecomposable tilting modules are classified
by their highest weights, see {\em loc.~cit.} Let $X_+$ be the set of
dominant weights, and for any $\lambda \in X_+$ let $T(\lambda)$
denote the indecomposable tilting module with highest weight $\lambda$. The
tensor product of tilting modules is again a tilting module.
 
Let us introduce the following
preorder relation $\le_T$ on $X_+$: $\lambda \le_T \mu$ iff $T(\lambda)$
is a direct summand of $T(\mu)\otimes \mbox{(some tilting module)}$.
We say that $\lambda \sim_T \mu$ if $\lambda \le_T\mu$ and $\mu \le_T\lambda$.
Obviously, $\sim_T$ is an equivalence relation on $X_+$. The equivalence
classes are called {\em weight cells}. The set of weight cells
has a natural order induced by $\le_T$. It was shown in \cite{O1}
that the partially ordered set of weight cells is isomorphic to the partially
ordered set of {\em two-sided cells} in the affine Weyl group $W$
associated with
$(Y,X,\ldots)$ ($W$ is a semidirect product of the finite Weyl group $W_f$ with
the dilated coroot lattice $pY$).

 Let $G$ and $\fg$ be the simply connected 
algebraic group and the Lie algebra
(both over ${\mathbb C}$) associated to $(Y,X,\ldots)$, and let $\cN$
be the {\em nilpotent cone} in $\fg$, i.e. the 
variety of $ad-$nilpotent elements.
It is well known that $\cN$ is a union of finitely many $G-$orbits called
{\em nilpotent orbits}.
 Using the theory of {\em support varieties} one defines {\em Humphreys map}
$\cH:$ \{ the set of weight cells\} $\to$ \{ the set of closed $G-$invariant
subsets of $\cN$\}, see \cite{O2}. The construction is as follows: it is
known that cohomology ring $\Ho^{\bullet}(u)$ of small quantum group 
$u\subset U$ is isomorphic to the ring of regular functions on $\cN$ (this is
a Theorem due to V.~Ginzburg and S.~Kumar, see \cite{GK}); now let $A$ be a 
weight cell and take any weight 
$\lambda \in A$, then $\Ext^\bullet (T(\lambda ), T(\lambda ))$ is naturally
a module over $\Ho^\bullet (u)$, so it can be considered as a coherent sheaf on
$\cN$, finally Humphreys map $\cH (A)$ is just support of this sheaf.
The Conjecture due to J.~Humphreys says
that the image of map $\cH$ consists of irreducible varieties, i.e. the 
closures of
nilpotent orbits; moreover J.~Humphreys conjectured that this map coincides
with {\em Lusztig's bijection} between the set of two-sided cells in the affine
Weyl group and the set of 
nilpotent orbits, see \cite{H} and \cite{O2}. In particular,
Humphreys map should preserve {\em Lusztig's $a-$function}; this function is
equal to half of codimension in $\cN$ of the nilpotent orbit and is defined
purely combinatorially on the set of two-sided cells, see \cite{LuC4}.
The aim of this note is to show that the Humphreys map does not decrease
the $a-$function: for a weight cell $A$ corresponding to a two-sided
cell $\uc$ in $W$ we have the inequality $\codim_{\cN}\cH(A)\ge
a(\uc)$. This inequality follows easily from the definition of $\cH$, Theorem
4.1 in \cite{O2} and the Main Theorem below:

{\bf Main Theorem.} {\em Let $A$ be a weight cell corresponding to a
two-sided cell $\uc$ in the affine Weyl group. Then there exists a weight
cell $B\le_TA$ and a regular weight $\lambda \in B$ such that $\dim T(\lambda)$
is not divisible by $p^{a(\uc)+1}$ provided $p$ is sufficiently large.}

{\bf Remark.} It follows from the Humphreys' Conjecture proved by 
Bezrukavnikov that in the Main Theorem $B=A$.

The proof of this Theorem is based on formulas for characters of
indecomposable tilting modules obtained by W.~Soergel in \cite{S1,S2}.
In what follows we will freely use notations and results from \cite{S1}
and \cite{O2}.

{\bf Warning.} The character formulas for tilting modules use certain
Kazhdan-Lusztig elements in the Iwahori-Hecke algebra of $W$, and in
modules thereof. The Iwahori-Hecke algebra is an algebra over
${\mathbb Z}[v,v^{-1}]$, and we will only need its specialization at $v=1$.
So all the notions related with it (e.g. Kazhdan-Lusztig bases) will be
understood in the specialization $v=1$.

\section{Proof of the Main Theorem.}
Let $\le$ denote the Bruhat order on the affine Weyl group $W$. For any
$x\in W$ let $(-1)^x$ denote sign of $x$, that is $(-1)^{l(x)}$ where
$l(x)$ is length of $x$.

\subsection{}
We may and will suppose that our
root system $R$ is irreducible. Let $S$ be the set of simple reflections in the
affine Weyl group $W$. For any $s\in S$ let $W_s$ be the parabolic subgroup
generated by $S-\{ s\}$. The subgroup $W_s$ is finite. There exists a unique
point $p\mu_s\in X\otimes_{\BZ}\BQ$ invariant under the $W_s-$action.
In general, $\mu_s\not \in X$, but the denominators of its coordinates
contain only bad primes for $R$.

In particular, let $s_a\in S$ be the 
unique affine reflection. Then $W_{s_a}=W_f$
is the finite Weyl group. There exists a 
natural projection $W\to W_f, x\mapsto \ox$. This projection embeds
all the subgroups $W_s$ into $W_f$.

Recall from \cite{S1} that the set $W^f$ of minimal length representatives of
cosets $W_f\setminus W$ is identified with the set of dominant alcoves.

\subsection{}
\label{1.3}
Recall that any two-sided cell of $W$ intersects nontrivially some
$W_s$, see \cite{LuC4}. In the group algebra of $W_s$ there are two remarkable
bases: the Kazhdan-Lusztig base $\tuH_w, w\in W_s$,
and the dual Kazhdan-Lusztig base $\uH_w, w\in W_s$ (notations from \cite{S1}).
Recall that $\uH_w=\sum_{x\le w}p_{x,w}x$ and $\tuH_w=\sum_{x\le w}p_{x,w}
(-1)^{xw}x$ where $p_{x,w}$ are the values at 1 of Kazhdan-Lusztig polynomials.

Let $V=X\otimes_{\BZ}\BR$ be the reflection representation of $W_f$. For any
$s\in S$ the restriction of $V$ to $\overline{W}_s$ is isomorphic
to the reflection representation $V_s$ of $W_s$.

We refer the reader to \cite{LuC1} for definition and properties of Lusztig's
$a-$function. This function is defined on the set of elements of a Coxeter
group and takes values in $\bN \cup \infty$. We will use following properties
of $a-$function:

(i) $a$-function is constant on any two-sided cell, see \cite{LuC1} 5.4.

(ii) Suppose that $w\in W_0\subset W$ where $W_0$ is a parabolic subgroup of
$W$. Then values of $a-$function of $w$ calculated with respect to Coxeter
groups $W_0$ and $W$ coincide, see \cite{LuC2} 1.9 (d).

(iii) 
Let $w\in W_s$. The element $\tuH_w$ acts trivially on $S^i(V_s)$ for 
$i<a(w)$,see \cite{Lucl}. The space $S^{a(w)}(V_s)$ contains exactly one 
irreducible component (special representation) such that elements 
$\tuH_{w'}, w'\sim_{LR}w$, act nontrivially on it, see {\em loc.~cit}. 
Moreover, these elements generate an action of the full
matrix algebra on this component, see \cite{Luch} Chapter 5. We will say that
this special representation corresponds to $w$. 

{\bf Convention.} The equivalence relation $\sim_{LR}$ depends on the ambient
group, e.g. if $w_1, w_2 \in W_s$ then $w_1\sim_{LR}w_2$ in $W$ does not
imply $w_1\sim_{LR}w_2$ in $W_s$. In what follows the equivalence relation
$\sim_{LR}$ is considered with respect to $W_s$. In spite of this we apply
the notation $\le_{LR}$ with respect to $W$. We hope that this does not cause
ambiguity in what follows.

\subsection{}
Let $\Delta(\lambda)=\prod_{\alpha \in R_+}\frac{\langle \lambda ,
\alpha^{\vee}\rangle}{\langle \rho ,\alpha^{\vee}\rangle}$ be the Weyl
polynomial. For any $w\in W_s$ and $y\in W_f$ let us consider the
following polynomial in $\lambda$ and $\mu$:
$$
\Delta(y,W_s,w,\mu ,\lambda )=\sum_{x\le w}p_{x,w}\Delta(\mu +y\ox y^{-1}\lambda).
$$
{\bf Lemma.} {\em The lowest degree term of $\Delta(y,W_s,w,\mu ,\lambda )$ in
$\mu$ has degree $\ge a(w)$.}

{\bf Proof.} It is well known that the polynomial $\Delta (\lambda )$ is
skew-symmetric with respect to the $W_f-$action: $\Delta (y\lambda )=(-1)^y
\Delta (\lambda )$. Hence
$$
\Delta(y,W_s,w,\mu ,\lambda )=\Delta(1,W_s,w,y^{-1}\mu ,y^{-1}\lambda ).
$$
So for the proof of the Lemma it is enough to consider the case $y=1$.
Using skew-symmetricity of $\Delta (\lambda )$ with respect to the 
$W_f-$action once again, we have:
$$
\leqno{(*)}\quad \Delta(1,W_s,w,\mu ,\lambda )=
\sum_{x\le w}p_{x,w}(-1)^x\Delta
(\ox^{-1}\mu+\lambda).
$$
The element $\sum_{x\le w}p_{x,w}(-1)^xx^{-1}$ equals $\tuH_{w^{-1}}$, see
e.g. \cite{S1}, proof of the Theorem 2.7. Since $a(w)=a(w^{-1})$ (see 
\cite{LuC1}), the result follows from \ref{1.3}.

\subsection{} \label{2.4}
Now we consider $\Delta (\mu +\lambda)$ as polynomial
in two variables $\mu ,\lambda \in V$. The action of the Weyl group $W_f$
on the space $S^\bullet (V\oplus V)$ of all polynomials in two variables
$\mu$ and $\lambda$ via the variable $\mu$ is well-defined and preserves
degrees of polinomials with respect to both $\mu$ and $\lambda$.

{\bf Lemma.} {\em Let $W_s$ act on the
polynomial $\Delta (\mu +\lambda)$
by the rule $x\Delta (\mu +\lambda)=\Delta(\ox \mu +\lambda)$. Then the
representation generated by the 
summand of degree $a(w)$ in $\mu$ contains the special representation
corresponding to $w^{-1}$.}

{\bf Proof.} Let $E_1\subset S^{a(w)}(V)$ be the 
special representation of $W_s$ corresponding to
$w^{-1}$. According to \cite{Lucl} sect.3, the $W_f-$representation $E$
generated by $E_1$ is irreducible, occurs with multiplicity 1 in the space
of polynomials of degree $a(E_1)=a(w)$ and does not occur in the spaces of
polynomials of lower degree. Moreover, $E$ lies in the space of harmonic
polynomials which is identified with the cohomology of the 
flag variety $\Ho^{2a(w)}(G/B)$.
Hence the Lemma is reduced to the following statement:

\subsubsection{}
{\bf Lemma.} (R.~Bezrukavnikov) {\em Let $W_f$ act on the 
polynomial $\Delta (\mu +\lambda)$ by the
rule $x\Delta (\mu +\lambda)=\Delta(x\mu +\lambda)$. Then the representation
generated by the summand of degree $i$ in $\mu$ contains any irreducible
constituent of $\Ho^{2i}(G/B)$.}

{\bf Proof.} We identify cohomology space $\Ho^{\bullet}(G/B\times G/B)$ with
the space of harmonic polynomials in two variables $\mu$ and $\lambda$
(with respect to the group $W\times W$). It is well known that diagonal class
is represented by $\Delta (\mu +\lambda )$. Using Poincar\'e duality we
identify $\Ho^{\bullet}(G/B\times G/B)$ with $\End (\Ho^{\bullet}(G/B))$ (this
identification is not $W\times W-$equivariant since fundamental class is
$W-$antiinvariant but it is $W-$equivariant with respect to the action of
first copy of $W$). Now any vector $v\in \Ho^{\bullet}(G/B)$ defines 
$W-$equivariant map $\End (\Ho^{\bullet}(G/B))\to \Ho^{\bullet}(G/B),\;
x\mapsto xv$ where $W$ acts on $\End (\Ho^{\bullet}(G/B))=\Ho^{\bullet}(G/B)\otimes
(\Ho^{\bullet}(G/B))^*$ via the first factor and under this map $1\mapsto v$.
The diagonal class $\Delta (\mu +\lambda )$ corresponds to $1\in 
\End (\Ho^{\bullet}(G/B))$ and the summand of degree $i$ in $\mu$ corresponds
to $1\in \End (\Ho^{2i}(G/B))$. The result follows.
 
\subsection{}
\label{1.6}
Let $N$ denote the degree of the polynomial $\Delta(\lambda)$.

{\bf Lemma.} {\em Let us fix $\mu$ such that $\langle \mu ,\alpha^{\vee}\rangle
\ne 0$ for any $\alpha \in R$. Then there exists $w'\in W_s, w'\sim_{LR}w$,
 such that the 
summand of $\Delta (y,W_s,w',\mu ,\lambda)$ of degree $N-a(w)$ in $\lambda$
is nontrivial.}

{\bf Proof.} We may and will assume that $y=1$. By the for\-mu\-lae $(*)$ we
have $\Delta (1,W_s,w_1,\mu ,\lambda)=\tuH_{w_1^{-1}}\Delta (\mu +\lambda )$ 
where $W_s$ acts on $\Delta (\mu +\lambda )$ via the variable $\mu$. 
Since elements $\tuH_{w_1^{-1}},\; w_1\sim_{LR}w$ generate  an action of the 
full matrix algebra on the special representation corresponding to $w^{-1}$ by
\ref{1.3} the Lemma 
\ref{2.4} show that the set of summands of degree $a(w)$ in $\mu$ 
of $\tuH_{w_1^{-1}}\Delta (\mu +\lambda )$ where $w_1$ runs through all $w_1
\sim_{LR}w$ contains a basis over the field of rational functions in $\lambda$
of special representation corresponding to $w^{-1}$. Evidently these summands
are exactly summands of $\Delta (1,W_s,w_1,\mu ,\lambda)$ of degree $N-a(w)$
in the variable $\lambda$. 
 
Our Lemma claims that
this set contains at least one nonzero element when we specialize $\mu$ to
any weight satisfying conditions of the Lemma. Consider ideal generated by
this set in a ring of polynomials in $\mu$ with coefficients which are
rational functions in $\lambda$. Evidently, the Lemma is a
consequence of the following statement: 

\subsubsection{}
{\bf Lemma.} {\em Let $U$ be an irreducible $W_f$-submodule
of $S^{\bullet}(V)$ not contained in $(S^+(V))^{W_f}$.
In other words, $U$ projects nontrivially to $S^\bullet(V)/(S^+(V))^{W_f}=
\Ho^{2\bullet}(G/B)$. Then the zero set of the ideal of $S^{\bullet}(V)$
generated by $U$ is contained in the union of hyperplanes $\langle \mu
,\alpha^{\vee}\rangle =0, \alpha \in R$.}

{\bf Proof.} Evidently, the ideal generated by $U$ is $W_f-$invariant. 
By Poincar\'e duality for any $0\ne v\in \Ho^i(G/B)$ there exists
$v'\in \Ho^{2N-i}(G/B)$ such that $vv'$ represents fundamental class of $G/B$.
Hence the ideal generated by $U$ contains an element $\omega \in S^N(V)$ which
projects nontrivially on $\Ho^{2N}(G/B)$. The alternation $\omega'=
\frac{1}{|W_f|}\sum_{w\in W_f}(-1)^ww(\omega )$ is also contained in our ideal
and projects nontrivially on $\Ho^{2N}(G/B)$. But $\omega'$ should be a
nonzero multiple of Weyl polynomial $\Delta (\lambda )$ since Weyl polynomial
is unique up to scalar $W-$antiinvariant in $S^N(V)$. The Lemma is proved.

\subsection{}
Let $\uc \subset W$ be a two-sided cell.
Choose $W_s$ such that $W_s\cap \uc \ne
\emptyset$ (this is possible by \cite{LuC4} Theorem 4.8(d)). Let us fix $w_1\in
W_s\cap \uc$. We choose $y\in W^f$ minimal with the property:

(**) For some $w\in W_s$ such that $w\sim_{LR}w_1$ the summand of
$\Delta (\oy ,W_s,w,yp\mu_s,\oy \lambda)$ of degree $N-a(w)$ in $\lambda$
is nonzero.

By Lemma \ref{1.6} such $y$ exists since there exists $y\in W^f$ such
that $y\mu_s$ lies strictly inside the dominant Weyl chamber.

In the following Lemma we use notations of \cite{S1}.

{\bf Lemma.} {\em Let $y\in W^f$ and $w\in W_s$ be as above. Then
the element $\uN_y\uH_w\in \cN$ is a sum of elements
$\uN_x, x\le_{LR}\uc$, with positive integral coefficients, and hence
 can be considered as the character of
tilting module in a regular block.}

{\bf Proof.} By the formulae in the end of Proposition 3.4 of \cite{S1} we
have:
$$
\uN_1\uH_x=\left \{
\begin{array}{cc}
\uN_x & \mbox{if}\ x\in W^f\\
0 & \mbox{if}\ x \not \in W^f.\\
\end{array}
\right.
$$
So, $\uN_y\uH_w=\uN_1\uH_y\uH_w$ and the Lemma follows from 
the definition of cells, together with the positivity
properties of multiplication in the Iwahori-Hecke algebra, 
see e.g. \cite{LuC1} \S 3.

\subsection{Proof of the main Theorem}
We can rewrite the element $\uN_y\uH_w$ as
$$
\uN_y\uH_{w}=N_1\sum_{y_1\in W^f,y_1\le y}n_{y_1,y}\sum_{x\le w}p_{x,w}H_{y_1x}.
$$
Let $\lambda_1$ be a regular weight from the fundamental alcove. The dimension
of the tilting module $T$ in the linkage class of $\lambda_1$ with character
given by $\uN_y\uH_w$ is equal to
$$
\sum_{y_1\in W^f,y_1\le y}n_{y_1,y}\sum_{x\le w}p_{x,w}\Delta(y_1x\cdot 
\lambda_1+\rho).
$$
Now let us write $\lambda_1=-\rho+p\mu_s+\lambda$. 

We have
$$
\mbox{dim}T=
\sum_{y_1\le y}n_{y_1,y}\sum_{x\le w}p_{x,w}\Delta 
(y_1p\mu_s+\overline{y_1x}\lambda)=
$$
$$
=\sum_{y_1\le y}n_{y_1,y}\Delta(\oy_1,W_s,w,y_1p\mu_s,\oy_1\lambda).
$$
According to (**), for some $w\sim_{LR}w_1$ the polynomial $\mbox{dim}T$
has nonvanishing summand of degree $N-a(\uc)$ in $\lambda$. Hence, for
$p\gg0$ it is possible to choose such $\lambda$ that this summand is not
divisible by $p$ and $\lambda_1$ lies in the lowest alcove.
The Main Theorem is proved.

{\large Acknowledgements.} This note is a re\-sult of con\-ver\-sa\-tions 
with many ma\-the\-ma\-ti\-cians. Es\-pe\-cially I wish to thank 
R.~Bezrukavnikov, M.~Finkelberg, J.~Humphreys,
J.~C.~Jantzen and G.~Rybnikov for their generous help and extremely useful
discussions. I am grateful to the referee for careful reading of the paper
and useful comments.

\end{document}